\documentclass[12pt]{article}

\usepackage{amssymb}

\newtheorem {theorem} {Theorem}
\newtheorem {definition} [theorem]{Definition}
\newtheorem {proposition} [theorem]{Proposition}

\newtheorem {lemma}  [theorem]{Lemma}

\newcommand{\bbox}{\ \hfill\rule[-1mm]{2mm}{3.2mm}}

\title{\sc Integrability of planar polynomial differential systems
 through linear differential equations.\thanks{The second and third
 authors are partially supported by a MCYT grant number BFM
2002-04236-C02-01. The second author is partially supported by
DURSI of Government of Catalonia's Acci\'o Integrada ACI2002-24.}}

\author{{\sc H. Giacomini$^{\ (1)}$, J. Gin\'e$^{\ (2)}$ and M.
Grau$^{\ (2)}$}}

\date{}

\begin{document}
\maketitle \vspace{-0.5cm}

\begin{center}
\noindent {\normalsize{$^{\ (1)}$ Laboratoire de Math\'ematiques
et Physique Th\'eorique. C.N.R.S. UMR 6083. \\ Facult\'e des
Sciences et Techniques. Universit\'e de Tours. \\ Parc de
Grandmont 37200 Tours, FRANCE.
\\ E-mail: {\tt giacomin@celfi.phys.univ-tours.fr}
\\ $ $ \\
$^{\ (2)}$ Departament de Matem\`atica. Universitat de Lleida. \\
Avda. Jaume II, 69. 25001 Lleida, SPAIN. \\ {\rm E--mails:} {\tt
gine@eup.udl.es} and {\tt mtgrau@matematica.udl.es} }}
\end{center}

\begin{abstract}
In this work, we consider rational ordinary differential
equations $ dy/dx = Q(x,y)/P(x,y),$ with
$Q(x,y)$ and $P(x,y)$ coprime polynomials with real coefficients.
We give a method to construct equations of this type for which a first integral can be expressed from two independent solutions of a second--order homogeneous linear
differential equation. This first integral is, in general, given by a non Liouvillian function.
\par
We show that all the known families of quadratic systems with an
irreducible invariant algebraic curve of arbitrarily high degree and without a rational first integral can be constructed by using this method. We also present a new example of this kind of families.
\par
We give an analogous method for constructing rational equations but by means of a linear differential equation of first order.
\end{abstract}


\section{Introduction}

This paper deals with rational ordinary differential equations such as
\begin{equation}
\frac{dy}{dx} \ = \ \frac{Q(x,y)}{P(x,y)}\, , \label{eqrac}
\end{equation}
where $Q(x,y)$ and $P(x,y)$ are coprime polynomials with real
coefficients. We associate to this rational equation a planar
polynomial differential system by introducing an independent
variable $t$ usually called {\em time}. Denoting by $\dot{ }=
d/dt$, we have
\begin{equation}
\dot{x}=P(x,y), \quad \dot{y}=Q(x,y), \label{eq0}
\end{equation}
where $(x,y) \in \mathbb{R}^2$. This system defines the vector field ${\mathcal{X}} = P(x,y) \frac{\partial}{\partial x} + Q(x,y) \frac{\partial}{\partial y}$ over $\mathbb{R}^2$ and, equivalently, the $1$-form $\omega = Q(x,y) dx - P(x,y) dy$. We indistinctively talk about equation (\ref{eqrac}) and system (\ref{eq0}). Let {\rm d} be the maximum degree of $P$ and $Q$. We say that system (\ref{eq0}) is of degree {\rm d}. When ${\rm d}=2$, we say that (\ref{eq0}) is a {\em quadratic system}.
\par
In order to simplify notation, we define $\mathbb{R}[x,y]$ as the ring of polynomials in two variables with real coefficients and $\mathbb{R}(x,y)$ as the field of rational functions in two variables with real coefficients, that is, the quotient field of the previous ring. Analogous definitions stand for $\mathbb{R}[x]$ and $\mathbb{R}(x)$.
\newline

We have an equation (\ref{eqrac}) defined in a certain class of functions, in this case, the rational functions with real coefficients $\mathbb{R}(x,y)$ and we consider the problem whether there is a first integral in another, possibly larger, class. For instance, as we will discuss later on, H. Poincar\'e stated the problem of determining when a system (\ref{eq0}) has a rational first integral. A ${\mathcal{C}}^k$ function $H: {\mathcal{U}} \to \mathbb{R}$ such that it is constant on each trajectory of (\ref{eq0}) and it is not locally constant is called a {\em first integral} of system (\ref{eq0}) of class $k$ and the equation $H(x,y)=c$ for a fixed $c \in \mathbb{R}$ gives a set of trajectories of the system, but in an implicit way. When $k \geq 1$, these conditions are equivalent to $\omega \wedge dH = 0$ and $H$ not locally constant. The problem of finding such a first integral and the functional class it must belong to is what we call the {\em integrability problem}.
\par
To find an integrating factor or an inverse integrating factor for
system (\ref{eq0}) is closely related to finding a first integral
for it. When considering the integrability problem we are also
addressed to study whether an (inverse) integrating factor belongs
to a certain given class of functions.
\begin{definition}
Let ${\mathcal{W}}$ be an open set of $\mathbb{R}^2$. A function
$\mu: {\mathcal{W}} \to \mathbb{R}$ of class
${\mathcal{C}}^k({\mathcal{W}})$, $k>1$, that satisfies the linear
partial differential equation
\begin{equation}
\omega \wedge d \mu = \mu \, d \omega, \label{eqfi}
\end{equation}
is called an {\em integrating factor} of system (\ref{eq0}) on
${\mathcal{W}}$.
\end{definition}
It has been shown than an easier function to find which also gives
additional properties for a differential system (\ref{eq0}) is the
inverse of an integrating factor, that is, $V=1/\mu$, which is called
{\em inverse integrating factor}.
\par
We note that $\{V=0\}$ is formed by orbits of system (\ref{eq0}).
The function $\mu =1/V$ defines on ${\mathcal{W}} \setminus
\{V=0\}$ an integrating factor of system (\ref{eq0}), which allows
the computation of a first integral of the system on
${\mathcal{W}} \setminus \{V=0\}$. The {\it first integral $H$
associated to the inverse integrating factor} $V$ can be computed
through the integral
$
H(x,y)=  \int \omega/V,
$
and the condition (\ref{eqfi}) for $\mu=1/V$ ensures that this line integral is well defined.
\par
The inverse integrating factors play an important role in two of
the most difficult open problems of qualitative theory of
planar polynomial vector fields, which are the center problem and $16^{th}$
Hilbert problem. In \cite{CGGLl}, it has been noticed that for
many polynomial differential systems with a center at the origin
there is always an inverse integrating factor $V$ globally defined in all $\mathbb{R}^2$, which is
usually a polynomial. However, the first integral for a polynomial
differential system with a center at the origin can be very
complicated.
\newline

We say that a function $f(x,y)$ is an {\em invariant} for a system (\ref{eq0}) if $\omega \wedge df = k f$ with $k(x,y)$ a polynomial of degree lower or equal than ${\rm d}-1$, where {\rm d} is the degree of the system. This polynomial
$k(x,y)$ is called the {\em cofactor} of $f(x,y)$. In
the previous equality and all along this paper we use the
convention of identifying the space of functions over
$\mathbb{R}^2$ and the space of $2$--forms over $\mathbb{R}^2$. In case $f(x,y)=0$ defines a curve in the real plane, this definition implies that $\omega \wedge df$ equals zero on the points such
that $f(x,y)=0$. In case $f(x,y)$ is a polynomial we say that $f(x,y)=0$ is an {\em invariant algebraic curve} for system (\ref{eq0}).
\par
Let us consider $f(x,y)=0$ an invariant algebraic curve for system (\ref{eq0}), we will always assume that $f(x,y)$ is an irreducible polynomial in $\mathbb{R}[x,y]$. Otherwise, it can be shown that each of its factors is an invariant algebraic curve for system (\ref{eq0}). We will denote by $n$ the degree of the polynomial $f(x,y)$.
\par
In \cite{Darboux}, G. Darboux gives a method for finding an
explicit first integral for a system (\ref{eq0}) in case that
${\rm d}({\rm d}+1)/2 + 1$ different irreducible invariant
algebraic curves are known, where {\rm d} is the degree of the
system. In this case, a first integral of the form
$H=f_1^{\lambda_1}f_2^{\lambda_2} \ldots f_r^{\lambda_r}, $ where
$f_i(x,y)=0$ is an invariant algebraic curve for system
(\ref{eq0}) and $\lambda_i \in \mathbb{C}$ not all of them null,
for $i=1,2,\ldots,r$, $r \in \mathbb{N}$, can be defined in the
open set $\mathbb{R}^2 \setminus \Sigma$, where $ \Sigma = \{
(x,y) \in \mathbb{R}^2 \mid (f_1 \cdot f_2 \cdot \dots \cdot
f_r)(x,y) = 0 \}.$ The functions of this type are called {\em
Darboux} functions. We remark that, particularly, if $\lambda_i
\in \mathbb{Z}$ , $\forall i=1,2,\ldots,r$, $H$ is a {\em rational
first integral} for system (\ref{eq0}). In this sense J. P.
Jouanoulou \cite{Jouanolou}, showed that if at least ${\rm d}({\rm
d}+1) + 2$ different irreducible invariant algebraic curves are
known, then there exists a rational first integral.
\par
The main fact used to prove Darboux's theorem (and Jouanoulou's
improvement) is that the cofactor corresponding to each invariant
algebraic curve is a polynomial of degree $\leq {\rm d}-1$.
Invariant functions can also be used in order to find a first
integral for the system. This observation permits a generalization
of Darboux's theory which is given in \cite{garciagine}, where, for instance, non-algebraic invariant curves with an algebraic cofactor for a polynomial system of degree $4$ are presented. In our work, we give other families of systems with such invariant curves.
\par
C. Christopher, in \cite{Christopher1}, studies the multiplicity
of an invariant algebraic curve and gives the definition for
exponential factor, which is a particular case of invariant for system (\ref{eq0}).
\begin{definition}
Let $h,g$ be two coprime polynomials. The function $e^{h/g}$ is
called an {\em exponential factor} for system (\ref{eq0}) if for
some polynomial $k$ of degree at most ${\rm d}-1$, where {\rm d} is the
degree of the system, the following relation is fulfilled: $ \omega
\wedge d \left( e^{h/g} \right) = k e^{h/g}.$ As before, we say
that $k(x,y)$ is the {\em cofactor} of the exponential factor
$e^{h/g}$.
\end{definition}
We note that an exponential factor $e^{h/g}$ does not define an
invariant curve, but the next proposition, proved in
\cite{Christopher1}, gives the relationship between both notions.
\begin{proposition} {\sc \cite{Christopher1}} \  If $F=e^{h/g}$ is an exponential factor and $g$ is not a constant, then $g=0$ is an invariant algebraic curve, and $h$ satisfies the equation $\omega \wedge dh = (h k_g + g k_F)$ where $k_g$ and $k_F$ are the cofactors of $g$ and $F$, respectively.
\end{proposition}
\par
All these previous results are closely related to a result of M.
F. Singer \cite{singer} which represents an important progress in
the resolution of the integrability problem when considering first
integrals for a system (\ref{eq0}) in the class of Liouville
functions. Roughly speaking, we can define a {\em Liouville
function} or a function which can be expressed by means of
quadratures, as a function constructed from rational functions
using composition, exponentiation, integration, and algebraic
functions. A precise definition of this class of functions is
given in \cite{singer}.
\begin{theorem} {\sc \cite{singer}} \  Let us consider the polynomial $1$--form $\omega = Q dx - P dy$ related to system (\ref{eq0}). System (\ref{eq0}) has a Liouville first integral if, and only if, $\omega$ has an inverse integrating factor of the form $V = \exp \int \xi$, where $\xi$ is a closed rational $1$--form. \label{thsinger}
\end{theorem}
We notice that the conditions on the function $V$ given in this Theorem can be restated as $d \omega= \xi \wedge \omega$ and $d \xi =0$.
\par
Taking into account Theorem \ref{thsinger}, C. Christopher
\cite{Christopher2} gives the following result, which precises the form of the inverse integrating factor.
\begin{theorem} {\sc \cite{Christopher2}} \  If system (\ref{eq0}) has an inverse integrating factor of the form $ \exp \int \xi $ with $d \xi =0$ and $\xi= \xi_1 dx + \xi_2 dy$ where $\xi_i$, $i=1,2$, are rational functions in $x$ and $y$, then there exists an inverse integrating factor of system (\ref{eq0}) of the form
\[ V = \exp \{D/E \} \prod C_i^{l_i}, \]
where $D$, $E$ and the $C_i$ are polynomials in $x$ and $y$, and $l_i \in \mathbb{C}$.
\label{thchris}
\end{theorem}
Theorem \ref{thchris} states that the search for Liouvillian first integrals can be reduced to the search of invariant algebraic curves and exponential factors.
\newline

In \cite{Poincare5}, H. Poincar\'e stated the following problem concerning the integration of an equation (\ref{eqrac}): {\em Give conditions on the polynomials $P$ and $Q$ to recognize when there exists a rational first integral.} As the same H. Poincar\'e noticed, a sufficient condition to solve this problem consists on finding an upper bound for the degree of the invariant algebraic curves for a given system (\ref{eq0}). From Darboux's result, it is known that for every polynomial vector field, there exists an upper bound for the possible degrees of irreducible invariant algebraic curves. However, it is a hard problem to explicitly determine such an upper bound. Some bounds have been given under certain conditions on the invariant curves, see D. Cerveau and A. Lins Neto work \cite{CerveauLinsNeto}, or on the local behavior of critical points, see M. Carnicer's work \cite{Carnicer}.
\par
In this sense, A. Lins Neto conjectured \cite{LinsNeto} that a
polynomial system (\ref{eq0}) of degree {\rm d} with an invariant
algebraic curve of degree high enough (where this bound only
depends on {\rm d}) would have a rational first integral. This
conjecture has been shown to be false by several counterexamples.
J. Moulin-Ollagnier \cite{Moulin} gives a family of quadratic
Lotka-Volterra systems, each with an invariant algebraic curve of
degree $2 \ell$, where $\ell$ is the parameter of the family,
without rational first integral. A simpler example is given by C.
Christopher and J. Llibre in \cite{ChLl}. In \cite{counter} a
family of quadratic systems with an invariant algebraic curve of
arbitrarily high degree without a Darboux first integral nor a Darboux inverse integrating factor is given. All these counterexamples exhibit a Liouvillian first integral.
\par
The natural conjecture at this step, also given by A. Lins Neto, see \cite{LinsNeto2}, after the counterexample of J. Moulin-Ollagnier appeared, is that a polynomial system (\ref{eq0}) of degree {\rm d} with an invariant algebraic curve of degree high enough (where this bound only depends on {\rm d}) has a Liouvillian first integral.
\newline

In this work we show a relationship between solutions of a class of systems (\ref{eq0}) and linear homogeneous ordinary
differential equations of order $2$ of the
form
\begin{equation}
A_2 (x) \, w''(x) + A_1(x) \, w'(x) + A_0(x) \, w(x) =0,
\label{eq1}
\end{equation}
where $x \in \mathbb{R}$, $w'(x)=dw(x)/dx$ and $w''(x)=dw'(x)/dx$.
We only consider equations (\ref{eq1}) where $A_i(x) \in
\mathbb{R}[x]$ for $0 \leq i \leq 2$ and $A_2(x) \not\equiv 0$.
\par
By means of a change of variable we rewrite an equation
(\ref{eq1}) as a polynomial differential system such that it has
an invariant related to $w(x)$. In case $w(x)$ is a
polynomial we get an invariant algebraic curve.
\par
Moreover, we give an explicit first integral for all the systems
built up by this method by means of two independent solutions of equation (\ref{eq1}).
\par
We give analogous results for a linear homogeneous ordinary differential equation of order $1$ such as
\begin{equation}
w'(x) + A(x) w(x) =0, \label{eq1l}
\end{equation}
where $x \in \mathbb{R}$, $w'(x)=dw(x)/dx$ and $A(x) \in
\mathbb{R}(x)$. All these results are given in Section
\ref{sect2}.
\newline

In Section \ref{sect3} we consider all the families of quadratic
systems with an algebraic curve of arbitrarily high degree known
until the moment of composition of this paper and we show that
they all belong to the construction explained in Section
\ref{sect2}. The families of quadratic systems with an algebraic
curve of arbitrarily high degree studied in this paper are the
ones appearing in \cite{algeblin, counter, ChLl, Moulin} and
one example more first appearing in this work. This new example consists on a biparametrical family of quadratic systems, which we give an explicit expression of a first integral for, such that when one of the parameters is a natural number, say $n$, the system exhibits an irreducible invariant algebraic curve of degree $n$.
\par
We give the explicit expression for the first integral of a
certain system (\ref{eq0}) by means of invariant functions for it,
and applying the Generalized Darboux's Theory as explained in
\cite{garciagine} where a new kind of first integrals, not only
the Liouvillian ones as in classical theories, is described. We
exemplify this result with the families of systems depending on
parameters described in Section \ref{sect3}. We remark that the
first integrals that we give in Section \ref{sect2} are not, in
general, of Liouvillian type. However, these first integrals are
Liouvillian at the values of parameters which correspond to the
systems with algebraic solutions. In the Subsection
\ref{subcenter}, we give an example of a $3$-parameter family of
quadratic systems with a center at the origin which can be
constructed by the method appearing in Section \ref{sect2} from an
equation (\ref{eq1l}).
\par
A question suggested by these examples is whether there are
polynomial systems which are not reversible nor Liouvillian
integrable which have a center and can be integrated by means of
Theorem \ref{th2}, see Section 2. The work \cite{BCLN} is related
to this question as it gives an example of an analytic system, not
polynomial, with a center which is not reversible nor Liouvillian
integrable. All the known families of polynomial vector fields
with a center at the origin are either Liouvillian integrable or
reversible, see \cite{Z1,Z2} for the definition of reversibility.
In \cite{Z1,Z2}, \.{Z}o\l\c{a}dek classifies all the reversible
cubic systems with a center. The reversible systems may have a
first integral not given by Liouville functions or no explicit
form of a first integral may be known. For instance the reversible
system $\dot{x}=-y+x^4$, $\dot{y}=x$ has a first integral composed
by Airy functions, see \cite{garciagine}, and no Liouvillian first
integral exists. The system $\dot{x}=-y^3 + x^2 y^2/2$,
$\dot{y}=x^3$ is an example given by Moussu, see \cite{Moussu},
which has a center at the origin since it is a monodromic and
reversible singular point and no explicit first integral is known
for this system.
\par
Since some examples of polynomial systems, which can be integrated
by the method described in Section 2, appear after a birrational
transformation, another suggested open question is if all the
polynomial systems with a center are birrationally equivalent to
one derived from Theorem \ref{th2} or from Theorem \ref{th20l}.

\section{Homogeneous linear differential equations of \\ order $\leq 2$ and planar polynomial systems. \label{sect2}}

Let us consider a homogeneous linear differential equation of order $2$:
\begin{equation}
A_2(x) \, w''(x)+ A_1(x) \, w'(x) + A_0(x) \, w(x) =0, \label{eq2}
\end{equation}
where $w'(x)=dw(x)/dx$, $w''(x)=dw'(x)/dx$, $A_i(x) \in
\mathbb{R}[x]$, $i=0,1,2$, and $A_2(x) \not\equiv 0$.
\par
\begin{theorem}
Given $g(x,y)= g_0(x,y)/g_1(x,y)$ with $g_i(x,y) \in
\mathbb{R}[x,y]$, $g_1(x,y) \not\equiv 0$ and $\partial g /
\partial y \not\equiv 0$, each nonzero solution $w(x)$ of
equation (\ref{eq2}) is related to a finite number of solutions
$y=y(x)$ of the rational equation
\begin{equation}
\frac{dy}{dx}= \displaystyle \frac{A_0(x)\, g_1^2+ A_1(x)\, g_1 \,
g_0 + A_2(x) \, g_0^2 + A_2(x) \left(g_1 \, \frac{\partial
g_0}{\partial x} - g_0 \, \frac{\partial g_1}{\partial x}  \right)
} {A_2(x) \left( g_0 \frac{\partial g_1}{\partial y}  - g_1
\frac{\partial g_0}{\partial y} \right) }, \label{eq3}
\end{equation}
by the functional change $dw/dx=g(x,y)w(x)$, which implicitly
defines $y$ as a function of $x$. \label{th2}
\end{theorem}
{\em Proof.} Let us consider equation (\ref{eq2}) and the
functional change $dw/dx=g(x,y)w(x)$ where $y=y(x)$, that is, $y$
is implicitly defined as a function of $x$. This change may also
be written as $w(x)=\exp ( \int_{x_0}^x g(s,y(s)) ds ),$ where
$x_0$ is any constant, and it is injective. We see that it is not
necessarily bijective unless the maximum degree of $g_1(x,y)$ and
$g_0(x,y)$ in the variable $y$ equals $1$. But it defines a finite
number of functions $y(x)$.
\par
By this functional change, equation (\ref{eq2}) becomes
\[ w(x)\left(A_0(x) + g \, A_1(x) + g^2 \, A_2(x) + A_2(x)
\frac{dy}{dx} \frac{\partial g}{\partial y} + A_2(x)
\frac{\partial g}{\partial x}\right)=0. \] We have that $w(x)$ is
a nonzero solution of (\ref{eq2}) so this equation is equivalent
to the ordinary differential equation of first order (\ref{eq3}).
Therefore, each non-zero solution $w(x)$ of equation (\ref{eq2})
corresponds to a finite number of solutions $y=y(x)$ of the planar
polynomial system (\ref{eq3}). \bbox

\begin{theorem} We consider the $1$-form related to equation (\ref{eq3})
\begin{eqnarray*}
\omega & = & \left( A_0(x) \, g_1^2+ A_1(x) \, g_1 \, g_0 + A_2(x)
\, g_0^2 + A_2(x) \left(g_1\frac{\partial g_0}{\partial x} - g_0
\frac{\partial g_1}{\partial x}  \right) \right) \, dx
\\ & & -
 A_2(x) \left( g_0 \frac{\partial g_1}{\partial y}  - g_1
\frac{\partial g_0}{\partial y} \right)  \, dy .
\end{eqnarray*} Let $w(x)$ be any nonzero solution of
equation (\ref{eq2}). Then the curve defined by $f(x,y)=0$, with
$f(x,y):=g_1(x,y) w'(x)-g_0(x,y)w(x)$ is invariant for system
(\ref{eq3}) and has the polynomial cofactor \[ \begin{array}{ll}
k(x,y) = & \displaystyle \left( A_0 (x) \, \frac{\partial
g_1}{\partial y} + A_1(x)\, \frac{\partial g_0}{\partial y}
\right) g_1 + A_2(x)\, g_0 \frac{\partial g_0}{\partial y} \vspace{0.2cm}  \\
& \displaystyle + A_2 (x) \left( \frac{\partial g_1}{\partial y}
\frac{\partial g_0}{\partial x} - \frac{\partial g_0}{\partial y}
\frac{\partial g_1}{\partial x} \right). \end{array} \]
\label{th21}
\end{theorem}
{\em Proof.} Let us consider $f(x,y)$ as defined above and let us
compute $\omega \wedge df$:
\begin{eqnarray*}
\begin{array}{lll}
\omega \wedge d f & = &  \displaystyle \left( \frac{\partial
g_1}{\partial x} w'(x) + g_1 w''(x) -  \frac{\partial
g_0}{\partial x} w(x) - g_0 w'(x) \right) A_2(x) \left( g_0
\frac{\partial g_1}{\partial y} \right.
\\  \vspace{0.2cm} & & \displaystyle \left. - g_1  \frac{\partial
g_0}{\partial y} \right)  + \left( \frac{\partial g_1}{\partial y}
w'(x) - \frac{\partial g_0}{\partial y} w(x) \right)
 \left[ A_0(x) \, g_1^2 + A_1 (x) \, g_1 \, g_0 \right.
 \\ \vspace{0.2cm} & &  \displaystyle \left. + A_2(x) \, g_0^2 + A_2(x) \left( g_1
\frac{\partial g_0}{\partial x}
  - g_0 \frac{\partial g_1}{\partial x} \right) \right]
\end{array}
\end{eqnarray*}
\begin{eqnarray*}
\begin{array}{lll}
 & = & \displaystyle g_1 \left( g_0  \frac{\partial g_1}{\partial y} -
g_1  \frac{\partial g_0}{\partial y} \right) A_2(x) \, w''(x) +
\\ \vspace{0.2cm} & & \displaystyle g_1 \left( A_1(x) \, g_0 + A_0(x) \, g_1 \right)
\left( \frac{\partial g_1}{\partial y} w'(x) - \frac{\partial
g_0}{\partial y} w(x) \right) +
\\ \vspace{0.2cm} & & \displaystyle \left(A_2(x) \left(\frac{\partial g_0}{\partial x}
\frac{\partial g_1}{\partial y} - \frac{\partial g_1}{\partial x}
\frac{\partial g_0}{\partial y} \right) + A_2(x) \, g_0
\frac{\partial g_0}{\partial y} \right) \left( g_1 w'(x)-g_0 w(x)
\right)
\end{array}
\end{eqnarray*}
Since $w(x)$ is a solution of (\ref{eq2}), we can substitute
$A_2(x) w''(x)$ by $-A_1(x) w'(x) - A_0(x) w(x)$. Therefore,
\begin{eqnarray*}
\begin{array}{lll}
\omega \wedge d f & = & \displaystyle  \left[ \left( A_0 (x) \,
\frac{\partial g_1}{\partial y} + A_1(x) \, \frac{\partial
g_0}{\partial y} \right) \, g_1 \, + \right. \vspace{0.2cm} \\  &
& \displaystyle \left. A_2(x) \, g_0 \, \frac{\partial
g_0}{\partial y} + A_2 (x) \left( \frac{\partial g_1}{\partial y}
\frac{\partial g_0}{\partial x} - \frac{\partial g_0}{\partial y}
\frac{\partial g_1}{\partial x} \right) \right] f(x,y).
\end{array}
\end{eqnarray*}
Then, we have that the function $f(x,y)$ is an invariant for system
(\ref{eq3}) and has the written polynomial cofactor.
\bbox

\begin{theorem} Let $\{ w_1(x), w_2(x) \}$ be a set of
fundamental solutions of equation (\ref{eq2}). We define
$f_i(x,y):=g_1(x,y) w_i'(x)-g_0(x,y)w_i(x)$, $i=1,2$. Then, system
(\ref{eq3}) has a first integral $H(x,y)$ defined by \[
H(x,y):=\frac{f_1(x,y)}{f_2(x,y)}=\frac{g_1(x,y) \, w_1'(x) \, -
\, g_0(x,y) \, w_1(x)}{g_1(x,y) \, w_2'(x) \, - \, g_0(x,y) \,
w_2(x)}.
\] \label{th22}
\end{theorem}
{\em Proof.} By Theorem \ref{th21}, we have that
$f_i(x,y):=g_1(x,y) w_i'(x)-g_0(x,y)w_i(x)$, $i=1,2$, are
invariants for system (\ref{eq3}), both with the polynomial
cofactor \[ k(x,y)=\left( A_0 (x) \, \frac{\partial g_1}{\partial
y} + A_1(x)\, \frac{\partial g_0}{\partial y} \right) g_1 +
A_2(x)\, g_0 \frac{\partial g_0}{\partial y} + A_2 (x) \left(
\frac{\partial g_1}{\partial y} \frac{\partial g_0}{\partial x} -
\frac{\partial g_0}{\partial y} \frac{\partial g_1}{\partial x}
\right). \] We remark that $f_1/f_2$ cannot be constant since the
two solutions $w_i(x)$, $i=1,2$, are independent. Therefore,
\[ \omega \wedge dH = \frac{f_2 (\omega \wedge df_1) - f_1 (\omega
\wedge df_2) }{f_2^2} = \frac{f_2 k f_1  - f_1 k f_2}{f_2^2} \equiv
0.
\] So, $H(x,y)$ is a first integral of system (\ref{eq3}). \bbox

\begin{lemma}
The function defined by \[ q(x):= A_2(x) \exp\left( \int_{x_0}^x
\frac{A_1(s)}{A_2(s)} ds \right) \] is an invariant for system
(\ref{eq3}), with cofactor $(A_1(x)+A_2'(x)) \left( g_0
\frac{\partial g_1}{\partial y}  - g_1 \frac{\partial
g_0}{\partial y} \right) $. \label{lem2}
\end{lemma}
We notice that $q(x)$ is a product of invariant algebraic curves
and exponential factors for system (\ref{eq2}), with complex
exponents.
\\
{\em Proof.} We compute $\omega \wedge d q$ and we have
\[ \omega \wedge d q =  \omega \wedge \frac{A_1(x) + A_2'(x)}{A_2(x)} \, q \, dx = (A_1(x)+A_2'(x)) \left( g_0
\frac{\partial g_1}{\partial y}  - g_1 \frac{\partial
g_0}{\partial y} \right)  \,  q .\] We notice that this algebraic
cofactor has degree $\leq {\rm d}-1$ provided that system
(\ref{eq3}) has degree {\rm d}. \bbox

\begin{proposition}
We use the same notation as in Theorem \ref{th2}. Let $w(x)$ be a
nonzero solution of (\ref{eq2}) and we define
$f(x,y):=w'(x)-g(x,y)w(x)$ and $q(x)$ as in Lemma \ref{lem2}. The
function $V(x,y)=q(x) f(x,y)^2$ is an inverse integrating factor
of system (\ref{eq3}). \label{prop2}
\end{proposition}
{\em Proof.} We only need to verify that $\omega \wedge dV + V d \omega =0$. We have that
\[
\begin{array}{ll} d \omega = & \displaystyle - \left[ 2 \, A_0(x) \, g_1 \, \frac{\partial
g_1}{\partial y} + A_1(x) \, \left( g_0 \frac{\partial
g_1}{\partial y} + g_1 \frac{\partial g_0}{\partial y} \right) + 2
\,
A_2 (x) \, g_0 \, \frac{\partial g_0}{\partial y} \right. \vspace{0.2cm} \\
& \left. \displaystyle  + 2 \, A_2(x) \left( \frac{\partial
g_1}{\partial y} \frac{\partial g_0}{\partial x}  - \frac{\partial
g_0}{\partial y} \frac{\partial g_1}{\partial x} \right) + A_2'(x)
\left( g_0 \frac{\partial g_1}{\partial y}   - g_1 \frac{\partial
g_0}{\partial y} \right) \right] .
\end{array}
\] Then,
\[
\begin{array}{rl}
\omega \wedge d V = & \omega \wedge (2 \, q \,f\, d f + f^2 \, d q
) =  2  \, q  \, f  \, (\omega \wedge d f) + f^2 \,  (\omega
\wedge d q) \vspace{0.2cm}
\\
 = & \displaystyle \left[ 2 \, A_0(x) \, g_1 \, \frac{\partial
g_1}{\partial y} + 2 \, A_1(x)
 \, g_1  \, \frac{\partial g_0}{\partial y} + 2  \, A_2(x) \, g_0 \, \frac{\partial g_0}{\partial y}
\right. \vspace{0.2cm}
\\
 & \displaystyle \left. + 2 \, A_2(x)
\left( \frac{\partial g_1}{\partial y} \frac{\partial
g_0}{\partial x}  - \frac{\partial g_0}{\partial y} \frac{\partial
g_1}{\partial x} \right) + A_1(x) \, g_0  \,  \frac{\partial
g_1}{\partial y} - A_1(x) \, g_1 \, \frac{\partial g_0}{\partial
y} \right. \vspace{0.2cm} \\ & \displaystyle \left.
 + A_2'(x) \left( g_0  \frac{\partial g_1}{\partial y} - g_1
\frac{\partial g_0}{\partial y} \right) \right] \, V
\vspace{0.2cm} \\ = & \displaystyle \left[ 2 \, A_0(x) \, g_1 \,
\frac{\partial g_1}{\partial y} + A_1(x) \, \left( g_0
\frac{\partial g_1}{\partial y} + g_1 \frac{\partial g_0}{\partial
y} \right) + 2 \, A_2 (x) \, g_0 \, \frac{\partial g_0}{\partial
y} \right. \vspace{0.2cm} \\ & \displaystyle \left. + 2  \, A_2(x)
\left( \frac{\partial g_1}{\partial y} \frac{\partial
g_0}{\partial x}  - \frac{\partial g_0}{\partial y} \frac{\partial
g_1}{\partial x} \right) + A_2'(x)  \left( g_0 \frac{\partial
g_1}{\partial y}   - g_1 \frac{\partial g_0}{\partial y} \right)
\right] \, V \vspace{0.2cm} \\ = & - V d \omega.
\end{array} \] \bbox

We remark that Theorem \ref{th22} gives, in general, a non
Liouvillian first integral for the planar polynomial systems
(\ref{eq3}). In Section \ref{sect3} we analyze some polynomial
systems constructed from Theorem \ref{th22} that have no
Liouvillian first integral.
\newline

We consider now a linear homogeneous ordinary differential equation of order $1$ such as
\begin{equation}
w'(x) + A(x) \, w(x) =0, \label{eq2l}
\end{equation}
where $x \in \mathbb{R}$, $w'(x)=dw(x)/dx$ and
$A(x)=A_0(x)/A_1(x)$ with $A_i(x) \in \mathbb{R}[x]$ and $A_1(x)
\not\equiv 0$. We give analogous results for this case whose
proofs are not given to avoid non useful repetitions.
\begin{theorem}
Given $g(x,y)= g_0(x,y)/g_1(x,y)$ with $g_i(x,y) \in
\mathbb{R}[x,y]$, $g_1(x,y) \not\equiv 0$ and $\partial g /
\partial y \not\equiv 0$ and $h(x)=h_0(x)/h_1(x)$ with
$h_i(x) \in \mathbb{R}[x]$ and $h_1(x) \not\equiv 0$, each nonzero
solution $w(x)$ of equation (\ref{eq2l}) is related to a finite
number of solutions $y=y(x)$ of the rational equation
\begin{equation}
\frac{dy}{dx} \, = \, \frac{ A_1(x) \,  h_0(x)  \,  g_1^2 - A_0(x)
\, h_1(x)  \, g_0  \, g_1 - A_1(x) \,  h_1(x) \left( g_1
\frac{\partial g_0}{\partial x} - g_0 \frac{\partial g_1}{\partial
x} \right) }{ A_1(x)  \,  h_1(x) \left( g_1  \frac{\partial
g_0}{\partial y} - g_0 \frac{\partial g_1}{\partial y} \right)} ,
\label{eq3l}
\end{equation}
by the functional change \[ w(x) = g(x,y) - \exp \left(- \int_0^x A(s) ds \right) \left[ \int_0^x \exp \left( \int_0^s A(r) dr \right) h(s) ds \right] . \]
\label{th20l}
\end{theorem}

\begin{theorem}
We consider the $1$-form related to equation (\ref{eq3l}) \[
\begin{array}{rl} \omega = & \displaystyle \left[ A_1(x)  \, h_0(x)  \, g_1^2 - A_0(x) \,  h_1(x)  \,  g_0  \, g_1
- A_1(x) \,  h_1(x) \left( g_1 \frac{\partial g_0}{\partial x} -
g_0
\frac{\partial g_1}{\partial x} \right) \right] \, dx \vspace{0.2cm} \\
& \displaystyle - \, A_1(x) \,  h_1(x) \left( g_1
\frac{\partial g_0}{\partial y} - g_0 \frac{\partial g_1}{\partial
y} \right)  dy .
\end{array} \] Let $w(x)$ be any nonzero solution of equation
(\ref{eq2l}), that is, for $C \in \mathbb{R}-\{ 0 \}$ we have
$w(x) = C \exp \left( - \int_0^x A(s) d s \right)$.
\par
Then, the function  \[ f(x,y) := g_1 \, w(x) - g_0 + g_1 \, \exp
\left( - \int_0^x A(s) ds \right) \left[ \int_0^x \exp \left(
\int_0^s A(r) dr \right) h(s) ds \right] \] is invariant for the
polynomial system (\ref{eq3l}), with the polynomial cofactor
\[
\begin{array}{ll}
k(x,y) = & \displaystyle -A_0(x) \, h_1(x) \, g_1 \,
\frac{\partial g_0}{\partial y} + A_1(x) \, h_0(x) \, g_1 \,
\frac{\partial g_1}{\partial y} \vspace{0.2cm} \\ & \displaystyle
+ A_1(x) \, h_1(x) \, \left( \frac{\partial g_0}{\partial y}
\frac{\partial g_1}{\partial x} - \frac{\partial g_1}{\partial y}
\frac{\partial g_0}{\partial x} \right).
\end{array}
\] \label{th21l}
\end{theorem}

\begin{lemma}
The function $q(x,y) = g_1(x,y) \exp \left( \int_0^x -A(s) ds
\right)$ is an invariant for system (\ref{eq3l}) with the same
polynomial cofactor as $f(x,y)$.  \label{lem2l}
\end{lemma}

\begin{theorem}
We use the same notation as in Theorem \ref{th21l} and Lemma
\ref{lem2l}. The function $H(x,y)$ defined by $H(x,y) : =
f(x,y)/q(x,y)$ is a first integral for system (\ref{eq3l}) and the
function $V(x,y) := A_1(x) \,  h_1(x)  \, g_1(x,y) \,  q(x,y)$ is
an inverse integrating factor. \label{th22l}
\end{theorem}
We remark that $H(x,y)$ is a Liouvillian function and, therefore, a system (\ref{eq3l}) has always a Liouville first integral.

In Section \ref{sect3} we give an example of a 3-parameter family
of quadratic systems with a center at the origin which can be
constructed following Theorem \ref{th20l}.

\section{Examples of families of quadratic systems.\label{sect3}}

\subsection{Quadratic systems with invariant algebraic curves of arbitrarily high degree linear in one variable.}

We first consider the examples of families of quadratic systems
with algebraic solutions of arbitrarily high degree appearing in
\cite{algeblin}. In that work all the invariant algebraic curves
linear in the variable $y$, that is, defined by $f(x,y) = p_1(x) y
+ p_2(x),$ where $p_1$ and $p_2$ are polynomials, are determined.
\par
The example appearing in \cite{counter} is a further study of an
example appearing in \cite{algeblin} and the example given in
\cite{ChLl} is also described in \cite{algeblin}. We show that all
these quadratic systems, with an invariant algebraic curve of
arbitrary degree can be constructed by the method explained in the
previous section. Moreover, we give the explicit expression of a
first integral for any value of the parameter $n$, even in the
case when $n$ is not a natural number. If $n$ is a natural number,
we obtain the invariant algebraic curves of arbitrary degree and a
Liouvillian first integral. However, when $n \not\in \mathbb{N}$
we obtain polynomial systems with a non Liouvillian first
integral.
\par
As it is shown in \cite{algeblin}, all these families of systems
can be written, after an affine change of variables if necessary,
in the form
\begin{equation}
 \dot{x}=\Omega_1(x), \quad
\dot{y}=(2 n+1)\, L'(x)\, \Omega_1(x)-\frac{n(n+1)}{2}\, \Omega_1(x) \, \Omega_1''(x)\,  - L(x)^2 + y^2, \label{eqalg0}
 \end{equation}
where $\Omega_1(x)$ is any quadratic polynomial, $L(x)$ is any
linear polynomial and $'=d/dx$. We have that system (\ref{eqalg0})
has an invariant curve $f(x,y)=0$, where $f(x,y):=p_1(x) y +
\Omega_1(x)p_1'(x)-L(x)p_1(x)$, with a cofactor $y+L(x)$, where
$p_1(x)$ is a solution of the second order linear differential
equation
\begin{equation}
\Omega_1(x) w''(x) + (\Omega_1'(x)-2 L(x)) w'(x) + \frac{n}{2} (4 L'(x)- (n+1) \Omega_1''(x)) w(x) = 0.
\label{eqalg1}
\end{equation}
In \cite{algeblin} it is shown that, in case $n \in \mathbb{N}$,
an irreducible polynomial of degree $n$ belonging to a family of
orthogonal polynomials is a solution of equation (\ref{eqalg1}).
For instance, when $\Omega_1(x)=1$, we get the Hermite
polynomials, when $\Omega_1(x)=x$, we get the Generalized Laguerre
polynomials and when $\Omega_1(x)=1-x^2$, we get the Jacobi
polynomials.
\newline

We consider again the general case in which $n \in \mathbb{R}$ and
we define $A_2(x):=\Omega_1(x)$, $A_1(x):=\Omega_1'(x)-2 L(x)$ and
$A_0(x):=\frac{n}{2} (4 L'(x)- (n+1) \Omega_1''(x))$. We have the
linear differential equation (\ref{eqalg1}) in the same notation
as in Theorem \ref{th2} and we consider $g(x,y):= \displaystyle
\frac{L(x)-y}{A_2(x)} $.
\newline

The system obtained by the method explained in Section \ref{sect2}
exactly coincides with system (\ref{eqalg0}). We consider a set of
fundamental solutions of equation (\ref{eqalg1}) $\{ w_1(x),
w_2(x) \}$ and applying Theorem \ref{th22}, we have a first
integral for system (\ref{eqalg0}) for any value of the parameter
$n \in \mathbb{R}$.
\par
In case $n \in \mathbb{N}$ we have that $w_1(x)$ degenerates to a polynomial and $w_1'(x)- g(x,y) w_1(x) = 0$ coincides with the algebraic curve given in the work \cite{algeblin}.
\newline

We explicitly give the first integral for each of the families
described in \cite{algeblin} and for $n \in \mathbb{R}$. We have
that $A_2(x)$ is a non-null quadratic polynomial in this case, and
depending on its number of roots, we can transform it by a real
affine change of variable to one of the following forms:
 $A_2(x) = 1, x, x^2, 1-x^2, 1+x^2$.
\newline

If $A_2(x)=1$, we can choose $L(x)=x$ by an affine change of coordinates. A set of fundamental
solutions $\{w_1(x), w_2(x) \}$ for (\ref{eqalg1}) with $n \in \mathbb{R}$ is
 \[
\begin{array}{ll}
w_{1}(x)= & \displaystyle 2^n\sqrt{\pi} \left(\frac{1}{\Gamma
\left(\frac{1-n}{2}\right)}\
{}_{1}F_1\left(-\frac{n}{2};\frac{1}{2};x^2 \right) -
 \frac{2x}{\Gamma\left(-\frac{n}{2}\right)}
 \ {}_{1}F_1
\left(\frac{1-n}{2};\frac{3}{2};x^2 \right) \right)\, ,
\vspace{0.2cm} \\
w_{2}(x)= & \displaystyle 2^n\sqrt{\pi} \left(\frac{1}{\Gamma
\left(\frac{1-n}{2}\right)}\
{}_{1}F_1\left(-\frac{n}{2};\frac{1}{2};x^2 \right) +
 \frac{2x}{\Gamma\left(-\frac{n}{2}\right)}
 \ {}_{1}F_1
\left(\frac{1-n}{2};\frac{3}{2};x^2 \right) \right)\, ,
\end{array}
\] where
$\Gamma(x)$ is the Euler's--Gamma function defined by $\Gamma(x)=
\int_{0}^{\infty} t^{x-1} e^{-t} dt$ and ${}_{1}F_1(a;b;x)$ is the
confluent hypergeometric function defined by the series
\[ {}_{1}F_1(a;b;x)= \displaystyle  \sum_{k=0}^{\infty}
\frac{(a)_k}{(b)_k} \frac{x^k}{k!},  \] with $(a)_k= a(a+1)(a+2)
\ldots (a+k-1)$, the Pochhammer symbol. See \cite{Abramowitz} for
further information about these functions. So, a first integral
for this system is the expression given in Theorem \ref{th22}:
$H(x,y):=f_1(x,y)/f_2(x,y)$, where
\[
\begin{array}{rl}
f_{1,2}(x,y) = & \displaystyle \pm \,
\Gamma\left(\frac{1-n}{2}\right) \left[   6 \, (x y-x^2+1) \
{}_1F_1\left(\frac{1-n}{2};\frac{3}{2};x^2\right) \right. \vspace{0.2cm} \\
& \displaystyle \left.
 -  \, 4 \, (n-1) \, x^2 \
{}_1F_1 \left(\frac{3-n}{2};\frac{5}{2};x^2\right)\right] +
 \vspace{0.2cm} \\ & \displaystyle
3  \, \Gamma\left(-\frac{n}{2}\right) \left[ 2nx \
{}_1F_1\left(1-\frac{n}{2};\frac{3}{2};x^2\right)+ (x-y)   \
{}_1F_1\left(-\frac{n}{2};\frac{1}{2};x^2\right)\right] .
\end{array}
\]
\par
When $n \in \mathbb{N}$, we have that (\ref{eqalg1}) corresponds to the equation for Hermite polynomials and $w_1(x)$ coincides with the Hermite polynomial of degree $n$. The invariant algebraic curve given in \cite{algeblin} corresponds to $f_1(x,y)=0$.
\newline

If $A_2(x)=x$, we choose $L(x)=\frac{1}{2}(x-\alpha)$, where $\alpha$ is an arbitrary real constant, and a set
of fundamental solutions for (\ref{eqalg1}) is:
 \[
 w_1(x) = \frac{(\alpha+1)_n}{\Gamma(n+1)} \,
{}_{1}F_1\left(-n;\alpha+1;x \right),\ \  w_2(x) = x^{-\alpha} \,
{}_{1}F_1\left(-\alpha-n;1-\alpha;x \right). \]
The first integral for this system is $H(x,y)= x^{\alpha} h_1(x,y) /h_2(x,y)$ with
\[
\begin{array}{ll}
h_1(x,y) = & \displaystyle
  (2 y - x + \alpha) \, (\alpha+1) \ {}_1F_1\left(-n;\alpha+1;x\right) - 2 \, n \, x \ {}_1F_1 \left(1-n;\alpha+2;x\right),
\vspace{0.2cm} \\ h_2(x,y) = & (2 y - x + \alpha)\, (\alpha - 1) \
{}_1F_1\left(-\alpha-n;1-\alpha;x\right) \vspace{0.2cm} \\ & -2\,
(\alpha+n)\, x \ {}_1F_1\left(1-\alpha-n;2-\alpha;x\right) .
\end{array}
\]
The first integral as given in Theorem \ref{th22} is
$f_1(x,y)/f_2(x,y)$ and we notice that $H(x,y)= c
f_1(x,y)/f_2(x,y)$ where $c \in \mathbb{R} - \{ 0 \}$. We do not
write $c$ in terms of the parameters of the system to simplify
notation.
\newline

When $n \in \mathbb{N}$, we have that (\ref{eqalg1}) corresponds
to the equation of Generalized Laguerre polynomials and $w_1(x)$
coincides with the Generalized Laguerre polynomial
$L_n^{(\alpha)}$. The invariant algebraic curve given in
\cite{algeblin} corresponds to $f_1(x,y)=0$, where
$f_1(x,y):=w_1'(x)-g(x,y)w_1(x)$.\\

If $A_2(x)=x^2$, the birrational transformation yet described in
\cite{algeblin}, $x=1/X$ and $y=(1/X)(1/2 - Y)$, makes this case
equivalent to the previous one.
\newline

If $A_2(x)=1-x^2$, we choose $L(x)=\frac{1}{2}((\alpha+\beta)x +
(\alpha-\beta))$, where $\alpha, \beta$ are two arbitrary real
constants, and a set of fundamental solutions for (\ref{eqalg1})
is:
 \[
\begin{array}{ccc} w_1(x)  & = & \displaystyle
\frac{(\alpha+1)_n}{\Gamma(n+1)} \
 {}_{2}F_1\left(-n,1+\alpha+\beta +
n;\alpha+1;\frac{1-x}{2}
 \right),
\end{array}
\]
\[
\begin{array}{ccc}
w_2(x) & = & \displaystyle (1-x)^{-\alpha} \
{}_{2}F_1\left(-\alpha-n,1+\beta+n;1-\alpha;\frac{1-x}{2} \right),
\end{array}
\]

where ${}_{2}F_1\left(a_1,a_2;b;x \right)$ is the hypergeometric
 function
defined by

\[ {}_{2}F_1\left(a_1,a_2;b;x\right)=  \displaystyle
 \sum_{k=0}^{\infty}
\frac{(a_1)_k (a_2)_k}{(b)_k} \frac{x^k}{k!}. \]

The first integral given in Theorem \ref{th22} is $H(x,y) =
(1-x)^{\alpha} h_1(x,y)/h_2(x,y)$, where
\[
\begin{array}{ll}
h_1 = & \displaystyle n \, (1+\alpha+\beta+n) \, (x^2- 1) \
{}_2F_1 \left( 1-n, 2+\alpha+\beta+n; 2+\alpha; \frac{1-x}{2}
\right)+ \vspace{0.2cm} \\ & \displaystyle
 \, (\alpha+1) \, ( (\alpha+ \beta)x+(\alpha-\beta) -2y)
\ {}_2F_1 \left(-n, 1+\alpha+\beta+n;1+\alpha;\frac{1-x}{2}
\right) , \vspace{0.2cm} \\ h_2 = & \displaystyle (\alpha-1) \,
((\alpha-\beta) x +(\alpha+b) + 2y) \ {}_2F_1 \left( -\alpha-n,
1+\beta+n;1-\alpha; \frac{1-x}{2} \right) + \vspace{0.2cm} \\ &
\displaystyle (\alpha+n)\, (1+\beta+n) \, (x^2-1) \
{}_2F_1\left(1-\alpha-n, 2+\beta+n;2-\alpha; \frac{1-x}{2} \right)
.
\end{array}\]
The first integral as given in Theorem \ref{th22} is
$f_1(x,y)/f_2(x,y)$ and we notice that $H(x,y)= c
f_1(x,y)/f_2(x,y)$ where $c \in \mathbb{R} - \{ 0 \}$. As before,
we do not write $c$ in terms of the parameters of the system to
simplify notation.
\par
When $n \in \mathbb{N}$, we have that (\ref{eqalg1}) corresponds
to the equation of Jacobi polynomials and  $w_1(x)$ coincides with
the Jacobi polynomial $P_n^{(\alpha,\beta)}(x)$ and the invariant
algebraic curve given in \cite{algeblin} corresponds to
$f_1(x,y)=0$, where $f_1(x,y):= w_1'(x) - g(x,y) w_1(x)$.\\

If $A_2(x)=1+x^2$ the complex affine change of variable $x={\rm
i}X$ makes this case equivalent to the previous one, as it is
shown in \cite{algeblin}. \newline

We have re-encountered by this method all the examples appearing
in \cite{algeblin} from a unified point of view. In addition, in
this work we have given an explicit expression of a first integral
for each case and for any value of the parameter $n \in
\mathbb{R}$. To this end, we have found invariants for the system
and we have applied the generalization of Darboux's method as
explained in \cite{garciagine} to be able to construct a first
integral which is, in general, of non Liouvillian type.

\subsection{A Lotka-Volterra system}

As it has been explained in the introduction, the first
counterexample to Lins Neto conjecture was given by J.
Moulin-Ollagnier in \cite{Moulin}. His example is a quadratic
system with two invariant straight lines and an irreducible
invariant algebraic curve $f(x,y)=0$ of degree $2 \ell$ when $\ell
\in \mathbb{N}$. This gives a family of systems depending on the
parameter $\ell$ which have a Darboux inverse integrating factor
when $\ell \in \mathbb{N}$ but no rational first integral. The
method used in \cite{Moulin} only shows the existence of such
invariant algebraic curve but no closed formula to compute it is
given. We give an explicit expression for an invariant by means of
Bessel functions for any value of $\ell \in \mathbb{R}- \{
\frac{1}{2} \}$ which, in the particular case $\ell \in
\mathbb{N}$ degenerates to the algebraic curve encountered in
\cite{Moulin}.
\par
We show that after a birrational transformation, this example
coincides with a system constructed with the method explained in
Section \ref{sect2}. A {\em birrational transformation} is a
rational change of variables whose inverse is also rational. This
kind of transformations bring polynomial systems such (\ref{eq0})
to polynomial systems and do not change the Liouvillian or non
Liouvillian character of the first integral.
\par
Let us consider the system appearing in \cite{Moulin} but assuming that $\ell \in \mathbb{R}- \{ \frac{1}{2} \}$
\begin{equation}
\dot{x}=x \left(1-\frac{x}{2} + y \right), \quad \dot{y}=y \left(-\frac{2 \ell +1}{2 \ell -1} + \frac{x}{2} - y \right).
\label{eqlv0}
\end{equation}
We make the birrational transformation \[ x= \frac{4 u v}{1 - 2 \ell}, \quad y= \frac{1 - 2 \ell}{4 v}, \] whose inverse is
\[ u = x y, \quad v= \frac{1 - 2 \ell}{4 y}. \]
By this transformation, system (\ref{eqlv0}) becomes
\begin{equation}
\dot{u}= \frac{2 u}{1 - 2 \ell}, \quad \dot{v}= \frac{1-2 \ell}{4} + \frac{2 \ell +1}{2 \ell -1} \, v + \frac{2 u}{2 \ell -1} \, v^2.
\label{eqlv1}
\end{equation}
We notice that the equation for the orbits satisfied by the variable $v$ as a function of $u$ is a Ricatti equation.
\par
Let us consider $g(u,v):=v$ and the linear differential equation of order $2$ given by
\begin{equation}
u \, w''(u) + \frac{1}{2}(1+2 \ell) \, w'(u)-\frac{1}{8} (1-2 \ell)^2 \, w(u) =0.
\label{eqlv2}
\end{equation}
Applying the method given in the previous section, this linear differential equation gives system (\ref{eqlv1}) modulus a change of time.
\par
A set of two fundamental solutions for equation (\ref{eqlv2}) is given by
\begin{equation}
w_1(u)= u^{(1-2\ell)/4} I_{\frac{1}{2}-\ell}\left( (1-2\ell)
\sqrt{\frac{u}{2}} \right), w_2(u)= u^{(1-2\ell)/4} I_{\ell
-\frac{1}{2}} \left((1-2\ell) \sqrt{\frac{u}{2}} \right),
\label{eqlv3}
\end{equation}
provided that $\ell$ is not of the form $\frac{1}{2}(1-2 r)$, with
$r$ an integer number, because in this case $w_1$ and $w_2$ are
linearly dependent. The function $I_\nu(u)$ is the Modified Bessel
function defined by the solution of the second order differential
equation
\begin{equation}
u^2 \, w''(u) + u \, w'(u) - (u^2+\nu^2) \, w(u) =0,
\label{eqlv4}
\end{equation}
and being bounded when $u \to 0$ in any bounded range of ${\rm arg} \, ( u)$ with $\mathbb{R}e(u) \geq 0$. See \cite{Abramowitz} for further information about this function.
\par
Hence, by Theorem \ref{th22} we have that $H(u,v)=f_1(u,v)/
f_2(u,v)$, where $f_i(u,v):=w_i'(u)-vw_i(u)$ for $i=1,2$, is a
first integral for system (\ref{eqlv1}). For a sake of simplicity
we consider the following renaming of the independent variable
$u=2 z^2/(1-2 \ell)^2$. This is not a birrational transformation
and that's why we only use it to simplify notation. The function
$H$ writes as:
\begin{equation}
H = \frac{ (1-2\ell)^2 I_{\left( \frac{1+2\ell}{2} \right)} (z) - 4 v z I_{\left( \frac{2\ell-1}{2} \right)} (z) }{ (1-2\ell)^2 I_{-\left( \frac{1+2\ell}{2} \right)} (z) - 4 v z I_{-\left( \frac{2\ell-1}{2} \right)} (z) }.
\label{eqlv5}
\end{equation}
\par
By Theorem \ref{th21} we have that $f_i(u,v)$, $i=1,2$ are
invariants with the same polynomial cofactor $k$ for system
(\ref{eqlv1}), so the curve $f(u,v)=0$ given by $f(u,v)=\pi z^{2
\ell+1} (f_1^2(u,v)- f_2^2(u,v))$ is also an invariant. We
multiply by $\pi$ only for esthetic reasons.
\par
Now we assume that $\ell \in \mathbb{N}$ and we show that $f=0$
defines an invariant {\em algebraic} curve. To this end we use the
following formulas appearing in \cite{Abramowitz,wolfram}. When
$\nu - \frac{1}{2} = n \in \mathbb{Z}$, we define $c(n)= -n\pi
\sqrt{-1}/2$ and the following relation is satisfied:
\begin{equation}
\begin{array}{ll}
I_\nu (z) = &  \displaystyle -\frac{1}{\sqrt{z}}\, {\rm e}^{c(n)}
\sqrt{\frac{2}{\pi}} \left\{ \sinh \left( c(n) - z \right)
\sum_{k=0}^{\lfloor \frac{2 | \nu | - 1}{4} \rfloor } \frac{(|\nu|
+ 2 k - \frac{1}{2}) ! }{(2 k)! (|\nu| - 2 k - \frac{1}{2} ) ! ( 2
z ) ^{2k}} \ \right. \vspace{0.2cm} \\ & \vspace{0.2cm}
\displaystyle \left. + \ \cosh \left( c(n) - z \right)
\sum_{k=0}^{\lfloor \frac{2 | \nu | - 3}{4} \rfloor} \frac{(|\nu|
+ 2 k + \frac{1}{2}) ! }{(2 k+1)! (|\nu| - 2 k - \frac{3}{2} ) ! (
2 z ) ^{2k+1}} \right\}\, ,
\end{array}
\label{eqlv6}
\end{equation}
where $\lfloor x \rfloor$ stands for the greatest integer $k$ such that $k \leq x$ and $| \nu |$ stands for the absolute value.
\par
From the former equation we obtain the following two equalities, with $\nu - \frac{1}{2} = n \in \mathbb{Z}$ and $\ell \in \mathbb{N}$,
\begin{equation}
\begin{array}{lll}
I_\nu^2(z) - I_{-\nu}^2(z) & = &  \displaystyle \frac{2}{\pi z} \sum_{k=0}^n (-1)^{k+1}
\frac{ (2n-k)! (2n - 2k)! }{k! ((n-k)!)^2 }
\left(\frac{1}{2z}\right)^{2(n-k)}\, ,
\end{array}
\label{eqlv7}
\end{equation}
\begin{equation}
\begin{array}{l}
I_{\ell + \frac{1}{2}}(z)  I_{\ell - \frac{1}{2}}(z) \  - \
I_{-(\ell + \frac{1}{2})}(z)\  I_{-(\ell - \frac{1}{2})}(z)  =
\displaystyle (-1)^{\ell} \frac{2}{\pi z} \vspace{0.2cm}   \\
\vspace{0.2cm} \displaystyle \  \left[ \left(  \displaystyle
\sum_{i=0}^{\lfloor \frac{\ell}{2} \rfloor} \frac{(\ell + 2
i)!}{(2i)! (\ell - 2 i)! } \left(\frac{1}{2z} \right)^{2i} \right)
\left( \displaystyle \sum_{j=0}^{\lfloor \frac{\ell - 2}{2}
\rfloor} \frac{(\ell + 2 j)!}{(2j+1)! (\ell - 2j -2)!}
\left(\frac{1}{2z} \right)^{2j-1} \right) - \right.  \\
\vspace{0.2cm}  \ \displaystyle \left. \left( \sum_{i=0}^{\lfloor
\frac{\ell - 1}{2} \rfloor} \frac{(\ell + 2 i + 1)!}{(2i+1)! (\ell
- 2i -1)!} \left(\frac{1}{2z} \right)^{2i-1} \right) \left(
\sum_{j=0}^{\lfloor \frac{\ell-1}{2} \rfloor} \frac{(\ell + 2 j
-1)!}{(2j)! (\ell - 2 j -1)! } \left(\frac{1}{2z} \right)^{2j}
\right) \right] .
\end{array}
\label{eqlv8}
\end{equation}
Then, we have that $f_1(z,v)= (1- 2\ell)^2\,
I_{\ell+\frac{1}{2}}(z) - 4 v z \, I_{\ell - \frac{1}{2}} (z) $
and $f_2(z,v)= (1-2 \ell)^2 \, I_{-(\ell+\frac{1}{2})}(z) - 4 v z
\,  I_{-(\ell - \frac{1}{2})} (z) $, and we write $f$ arranged in
powers of $v$:
\begin{eqnarray*}
f(z,v) & = &  \pi z^{2 \ell + 1} \left( (1-2\ell)^4 (I_{\ell+\frac{1}{2}} ^2 (z) - I_{-(\ell+\frac{1}{2})} ^2 (z) ) \right.
\\ & & \left.
- \ 8 v z (1- 2 \ell)^2 (I_{\ell + \frac{1}{2}}(z) I_{\ell - \frac{1}{2}}(z) - I_{-(\ell + \frac{1}{2})}(z) I_{-(\ell - \frac{1}{2})}(z) ) \right.
\\ & & \left.
+ 16 v^2 z^2 (I_{\ell-\frac{1}{2}} ^2 (z) - I_{-(\ell-\frac{1}{2})} ^2 (z) ) \right).
\end{eqnarray*}
Let us consider each coefficient of $v$ in $f(z,v)$ separately and we will show that it is an even polynomial in the variable $z$. The coefficient in $f(z,v)$ of $v^0$ is:
\[ \pi z^{2 \ell + 1} (1-2\ell)^4 (I_{\ell+\frac{1}{2}}^2 (z) - I_{-(\ell+\frac{1}{2})}^2 (z) ), \]
which by equation (\ref{eqlv7}) is an even polynomial in the variable $z$ of degree $2 \ell$. The coefficient in $f(z,v)$ of $v^2$ is:
\[ 16 \pi z^{2\ell + 3} (I_{\ell-\frac{1}{2}}^2 (z) - I_{-(\ell-\frac{1}{2})}^2 (z) ), \]
which also by equation (\ref{eqlv7}) is an even polynomial in the
variable $z$ of degree $2 \ell + 2$. Finally, the coefficient in
$f(z,v)$ of $v^1$ is:
\begin{equation}
8 \pi (1 - 2 \ell)^2 z^{2 \ell +2} (I_{\ell + \frac{1}{2}}(z)
I_{\ell - \frac{1}{2}}(z) - I_{-(\ell + \frac{1}{2})}(z) I_{-(\ell
- \frac{1}{2})}(z) ),
\end{equation}
which by equation (\ref{eqlv8}) is an even polynomial in the
variable $z$ of degree $2 \ell$.
\par
Hence, we have that $f(z,v)$ is an even polynomial in the variable
$z$ of total degree $2 \ell +4$. When rewriting $z= (1 - 2 \ell)
\sqrt{u}/ \sqrt{2}$ we have that $f(u,v)$ is a polynomial of total
degree $\ell+2$ which is irreducible. The fact of being
irreducible is easily seen because it is a polynomial of degree
two in $v$ and it cannot be decomposed in linear factors (the
discriminant is not a polynomial raised to an even power) and the
coefficients of $v^0$ and $v^2$ do not have any root in common.
\par
Undoing the birrational change of variables we get that $f(x,y)$ is an irreducible
polynomial of degree $2 \ell$ given by:
\begin{eqnarray*}
f(x,y) & = & x^{\ell + \frac{1}{2}} y^{\ell - \frac{1}{2}} \left[
2y \left( I_{\ell+\frac{1}{2}}^2 (z) -I_{-(\ell+\frac{1}{2})}^2
(z) \right) + \right.
\\  & & \left. 2 \sqrt{2} \sqrt{x y} \left( I_{\ell+\frac{1}{2}}
(z) I_{\ell-\frac{1}{2}} (z) - I_{-(\ell+\frac{1}{2})}
(z) I_{-(\ell-\frac{1}{2})} (z) \right) + \right. \\
& & \left. x \left( I_{\ell-\frac{1}{2}}^2 (z) -
I_{-(\ell-\frac{1}{2})}^2 (z) \right) \right],
\end{eqnarray*}
where $z$ is the same variable as before, that is, $z= \frac{(1-2\ell) \sqrt{x y}}{\sqrt{2}}$.
\par
By equation (\ref{eqlv5}) we can write the first integral for system (\ref{eqlv0}) for any value of $\ell \in \mathbb{R} - \{ \frac{1}{2} (1-2 r) \, \mid \, r \in \mathbb{N} \}$:
\[
H(x,y) = \displaystyle \frac{ \sqrt{2 y} \, I_{\left(
\frac{1+2\ell}{2} \right)} (z) - \sqrt{x} \, I_{\left(
\frac{2\ell-1}{2} \right)} (z) }{ \sqrt{2 y} \, I_{-\left(
\frac{1+2\ell}{2} \right)} (z) - \sqrt{x} \, I_{-\left(
\frac{2\ell-1}{2} \right)} (z) }.
\]

We have studied system (\ref{eqlv0}) for any value of the
parameter $\ell \in \mathbb{R} - \{ \frac{1}{2} (1-2 r) \, \mid \,
r \in \mathbb{N} \}$ giving an explicit expression for a first
integral using Theorem \ref{th2} and the Generalized Darboux's
theory as explained in \cite{garciagine}. This first integral is
not of Liouvillian type. Moreover, we give one of its invariants
with a polynomial cofactor. In the particular case $\ell \in
\mathbb{N}$, this invariant is the invariant algebraic curve whose
existence was proved in \cite{Moulin}.

\subsection{A new example of a family of quadratic systems with an invariant algebraic curve of arbitrarily high degree}

We give another example of a family of quadratic systems with an irreducible invariant algebraic curve of degree $2\ell$ when $\ell \in \mathbb{N}$, where $\ell$ is a parameter of the family. This family also depends on the parameter $a \in \mathbb{R}$.
\par
Let us consider the quadratic system
\begin{equation}
\begin{array}{lll}
\dot{x} & = & (2 a - 1) \ell x - a (2 \ell -1 ) y + 2 a (a - \ell) (2 \ell - 1) x^2 - 2 a^2 (2 \ell - 1)^2 x y,
 \vspace{0.2cm} \\ \vspace{0.2cm}
\dot{y} & = & y ( 2 (2 a -1) \ell + 2 a ( 2 a - 2 \ell -1) (2 \ell
-1) x - 4 a ^2 ( 2 \ell - 1)^2 y),
\end{array}
\label{eqnou1}
\end{equation}
where $a, \ell \in \mathbb{R}$ which satisfy $a \neq 0$, $\ell \neq \frac{1}{2}$ and $(2 \ell - 1) a^2 - 2 \ell \neq 0$.  An straightforward computation shows that system (\ref{eqnou1}) has $y=0$ and $y-x^2=0$ as invariant algebraic curves.
\par
Let us consider the following birrational transformation $x=Y$, $y=X Y^2$ whose inverse is $X= y/x^2$ and $Y=x$. In these new variables system (\ref{eqnou1}) becomes
\begin{equation}
\begin{array}{lll}
\dot{X} & = & 2 a (2 \ell -1) (X-1) X Y,  \vspace{0.2cm} \\
\dot{Y} & = & ((2 a -1) \ell + a ( 2 \ell -1) ( 2 a - 2 \ell - X)
Y - 2 a^2 ( 2 \ell -1)^2 X Y^2) \, Y.
\end{array}
\label{eqnou2}
\end{equation}
By a change of the time variable we can divide this system by $Y$
and the resulting system coincides with the one described in
Theorem \ref{th2} taking $A_2(X) := 2 X (X-1)^2$, $A_1(X) := ( 2 \ell
- 2 a + 3 X) (X-1)$, $A_0(X) := \ell (1 - 2 a)$ and $g(X,Y) := a (2
\ell -1) Y/ (X-1)$. The equation $A_2(X) w''(X) + A_1(X) w'(X) +
A_0(X) w (X) =0$ has the following set of fundamental solutions in
this case:
\begin{eqnarray*}
w_1(X) & = & (X-1)^{- \ell} \ {}_{2}F_1 \left( \frac{1}{2} - \ell, -\ell; a-\ell; X \right), \\
w_2 (X) & = & (X-1)^{- \ell} X^{1-a+\ell} \ {}_{2}F_1 \left( 1-a,
\frac{3}{2}-\ell; 2-a+\ell; X \right).
\end{eqnarray*}
By Theorem \ref{th21}, $f_i(X,Y)= w_i'(X) - g(X,Y) w_i(X)$,
$i=1,2$, define invariants with a polynomial cofactor for system
(\ref{eqnou2}). Moreover, by Theorem \ref{th22} we have a non
Liouvillian first integral given by $H(X,Y) = f_1(X,Y) /
f_2(X,Y)$.
\newline

In the particular case $\ell \in \mathbb{N}$, we notice that
$f_1(X,Y)=0$ is a rational function. It is an easy computation to
show that this rational function is a polynomial when rewritten in
coordinates $x$ and $y$. This polynomial gives place to an
invariant algebraic curve of degree $2 \ell$ for system
(\ref{eqnou1}). That is, by undoing the birrational
transformation, we deduce that $f_1(x,y)$ is an irreducible
invariant algebraic curve for system (\ref{eqnou1}), given by:
\begin{eqnarray*}
f_1(x,y) & = & 2(a- \ell) (\ell + (2 \ell -1) a x) x^{2 \ell - 1}
\ {}_2F_1 \left( \frac{1}{2} -\ell, - \ell; a - \ell;
\frac{y}{x^2} \right) \\ & & + \ell ( 2 \ell - 1) x^{2 \ell - 3}
(x^2 - y)\ {}_2F_1 \left( \frac{3}{2} - \ell, 1-\ell; 1 + a -\ell;
\frac{y}{x^2} \right) .
\end{eqnarray*}
It is easy to see that the polynomial  $f_1(x,y)$ has degree $2
\ell$ and the cofactor associated to the invariant algebraic curve
$f_1(x,y)=0$ is $\ell (2 \ell - 1) ((2a-1) + 4 a (a -\ell) x - 4
(2 \ell - 1) a^2 y)$.
\par
The first integral for (\ref{eqnou1}) is given by $H(x,y)=
y^{a-\ell} f_1(x,y)/ h(x,y)$, where
\[
\begin{array}{rl}
h(x,y) = & \displaystyle 2 \, (a-\ell -2) \, \left[ (a -\ell -1)
\, x^2 \,+ \right. \vspace{0.2cm} \\ & \displaystyle \left.
 (1-a-a
(2 \ell -1) x ) \, y \right] \, x^{7-2 a} \ {}_2F_1 \left( 1-a,
\frac{3}{2}-a; 2-a-\ell;\frac{y}{x^2} \right) +
 \vspace{0.2cm} \\ & \displaystyle
(a-1) \, (2 a -3) \, x^{5-2a} \, (x^2-y) \, y \
{}_2F_1\left(2-a,\frac{5}{2} - a;3-a+\ell; \frac{y}{x^2} \right).
\end{array}
\]

We notice that when both $a$ and $\ell$ belong to the set of
natural numbers, we have that $h(x,y)=0$ is an invariant algebraic
curve different from $f_1(x,y)=0$. Then we have a quadratic system
with a rational first integral $H(x,y)$ with arbitrary degree.

\subsection{A complete family of quadratic systems with a center at the origin \label{subcenter}}

In this subsection we give an example of a $3$-parameter family of
quadratic systems with a center at the origin which can be
constructed using Theorem \ref{th20l}. The family encountered corresponds
to the reversible case, see \cite{Schlomiuk1}.
\par
The family of quadratic systems depends on twelve parameters, but
up to affine transformations and positive time rescaling, we get a
family of five essential parameters. We have taken a system
(\ref{eq3l}) and we have chosen $g(x,y):=y^2$, $h(x) := 2 x (d \,
x - 1)/ (1 + a x)$ and $A(x):=2b/(1+ax)$, where $a,b,d$ are real
parameters. Using Theorem \ref{th20l}, we have encountered the
$3$-parameter family of quadratic systems next described. We
remark that in spite of the simplicity of the chosen polynomials
$g(x,y)$, $A(x)$ and $h(x)$, we amazingly obtain the complete
family of quadratic systems with a reversible center at the
origin. We notice that other choices of the functions $g(x,y)$,
$A(x)$ and $h(x)$ would give place to other families of polynomial
systems.
\par
Let us recall that a {\em center} is an isolated singular point of
an equation (\ref{eq0}) with a neighborhood foliated of periodic
orbits. When the linear approximation of an equation (\ref{eq0})
near a singular point has non-null purely imaginary eigenvalues,
the point can be a center or a focus. To distinguish between these
two possibilities is the so-called {\em center problem}. H.
Poincar\'e gave a method to solve it by defining a numerable set
of values, called Liapunov-Poincar\'e constants, which are all
zero when the singular point is a center and at least one of them
is not null when it is a focus. When these constants are computed
from a family of systems, they are polynomials on the coefficients
of the family. Hilbert's Nullstellensatz ensures that there always
exists a finite number of independent polynomials which generates
the whole ideal made up with all these Liapunov-Poincar\'e
polynomials. The zero-set of these independent polynomials gives
place to the center subfamilies. The reader is referred to
\cite{Schlomiuk1, Schlomiuk2} for a survey on this subject.
\par
The computation of these center cases for the family of quadratic
systems was done by Dulac \cite{Dulac} for the case of complex
systems and a proof for real systems is given in
\cite{lunkevichsibirskii}. We also refer Bautin \cite{Bautin} who
showed the existence of only three independent constants. The
computation of the zero set of these three independent values
gives place to four complete families of quadratic systems with a
center at the origin which are described in \cite{Schlomiuk2}.
\newline

Let us now consider an equation (\ref{eq3l}) such as $ (1 + a x)
\, w'(x) + 2 b \, w(x) = 0$ and $g(x,y)$ and $h(x)$ as formerly
defined. The rational equation as constructed in Theorem
\ref{th20l} is
\[ \frac{dy}{dx} = \frac{ - x + d x^2 - b y^2}{y + a x y} , \]
which gives the corresponding quadratic planar system
\begin{equation}
\dot{x}= y + a x y, \quad \dot{y}= - x + d x^2 - b y^2.
\label{eq20}
\end{equation}
We suppose that $a b (a+b) (a + 2 b) (a + b + d) \neq 0$. In case this value is zero, the origin of system (\ref{eq20}) is still a center but with a Darboux integrating factor instead of a Darboux first integral. This particular case can also be studied by our method, but we do not write it to avoid giving examples without essential differences.
\par
By Theorem \ref{th21l} we have that $f(x,y)$ is an invariant of system (\ref{eq20}) with cofactor $- 2 b y $, where $f(x,y)$ is given by
\begin{equation}
\begin{array}{lll}
f(x,y) & = & \vspace{0.2cm}  \displaystyle  b \, (a+b) \, (a + 2 b) \,
w(x) - (a +b +d) \, \displaystyle (1 + a x)^{-\frac{2b}{a}}  \\
\vspace{0.2cm} & &
 - b \, (a+b) \, (a + 2 b) \, y^2 \, + \, b \, (a+ 2 b) \, d  \, x^2
\\ \vspace{0.2cm} & &  -  2\, b  \, (a + b + d) \, x \, + \, a \, + \, b \, + \, d \, ,

\end{array}
\label{eq21}
\end{equation}
with $w(x)$ any non-zero solution of $ (1 + a x) w'(x) + 2 b w(x) = 0$, that is, $w(x)= C (1 + a x)^{-\frac{2b}{a}}$.
\par
Choosing $C= (a + b + d) \, (b (a+b) (a + 2 b) )^{-1}$ we get an
invariant conic. System (\ref{eq20}) has two invariant algebraic
curves, the former conic with cofactor $-2 b y$ and an invariant
straight line given by $1+a x=0$ with cofactor $y$. The Darboux
first integral
\[ H(x,y) = (1 + a x)^{\frac{2 b}{a}} f(x,y) \] coincides with the
first integral described in Theorem \ref{th22l}.
\par
The origin of this system is a center since it is a monodromic
singular point with a continuous first integral defined in a
neighborhood of it. This example addresses to the thought that
other families of polynomial systems of higher degree with a
center at the origin can be easily obtained by this method,
avoiding the cumbersome computation of Poincar\'e-Liapunov
constants.


\begin{thebibliography}{99}

 \bibitem{Abramowitz} {\sc Abramowitz, M.; Stegun, I.A.},
 {\it  \ Handbook of mathematical functions with formulas, graphs,
 and mathematical tables.} Edited by Milton Abramowitz and Irene A.
 Stegun. Reprint of the 1972 edition. Dover Publications, Inc.,
New York, 1992.

\bibitem{BCLN} {\sc Berthier, M.; Cerveau, D.; Lins Neto, A.}, {\it \ Sur les feuilletages analytiques r\'eels et le probl\`eme du
centre}, J. Differential Equations {\bf 131} (1996), 244--266.

\bibitem{Bautin} {\sc Bautin, N.N.},{\it \ On the number of limit
cycles which appear with the variation of coefficients from an
equilibrium position of focus or centre type.} Trans. Amer. Math.
Soc. {\bf 100} (1954), 397--413.

\bibitem{Carnicer} {\sc Carnicer, M.}, {\it \ The Poincar\'e problem in the nondicritical case}, Ann. Math. {\bf 140} (1994), 289--294.

\bibitem{CerveauLinsNeto} {\sc Cerveau, D.; Lins Neto, A.},{\it \ Holomorphic foliations in ${\mathbb{CP}}(2)$ having an invariant algebraic curve}, Ann. Inst. Fourier {\bf 41} (1991), 883--903.

 \bibitem{CGGLl} {\sc Chavarriga, J.; Giacomini, H.;
Gin\'e, J.; Llibre, J.},{\it \ On the integrability of
two-dimensional flows.} J. Differential Equations {\bf 157}
(1999), 163--182.

\bibitem{algeblin} {\sc Chavarriga, J.; Grau, M.},{\it \ Invariant algebraic curves linear in one variable for planar real quadratic systems.}, Appl. Math. Comput. {\bf 138} (2003), 291--308.

 \bibitem{counter} {\sc
Chavarriga, J.; Grau, M.}, {\it \ A family of quadratic polynomial
differential systems non-integrables Darboux and with algebraic
solutions of arbitrarily high degree}, submitted to Appl. Math.
Lett.

 \bibitem{Christopher1} {\sc Christopher, C.},{\it
\ Invariant
 algebraic curves and conditions for a centre.} Proc. Roy. Soc.
Edinburgh Sect. A {\bf 124} (1994), 1209--1229.

\bibitem{Christopher2} {\sc Christopher, C.},{\it \ Liouvillian first
integrals of second order polynomial differential systems.}
 Electronic
Journal of Differential Equations, Vol. {\bf
 1999}(1999), 1--7.

\bibitem{ChLl} {\sc Christopher, C.; Llibre, J.},{\it \ A family
 of quadratic
polynomial differential systems with invariant
 algebraic curves of
arbitrarily high degree without rational first
 integrals.} Proc. Amer. Math.
Soc. {\bf 130} (2002),  2025--2030 (electronic).

\bibitem{Darboux}{\sc Darboux, G.},{\it \ M\'emoire sur les
 \'equations
diff\'erentielles alg\'ebriques du premier ordre et du
 premier degr\'e
(M\'elanges)}, Bull. Sci. Math. {\bf 32} (1878),
 60--96; 123--144; 151--200.

\bibitem{Dulac} {\sc Dulac, H.},{\it \ D\'etermination et
int\'egration d'une certaine classe d'\'equations
diff\'erentielles ayant pour point singulier un centre,} Bull.
Sciences Math. S\'er. (2) {\bf 32}(1) (1908), 230--252.

\bibitem{garciagine} {\sc Garc\'{\i}a, I.A.; Gin\'e, J.},{\it \
 Generalized
cofactors and nonlinear superposition principles,} Appl. Math. Lett., to appear.

\bibitem{Jouanolou} {\sc Jouanolou, J.P.},{\it \
\'Equations de
 Pfaff alg\'ebriques.} Lecture Notes in Mathematics, {\bf 708}.
Springer, Berlin, 1979.

\bibitem{LinsNeto} {\sc Lins Neto, A.},{\it \ Folhea\c{c}\={o}es
em $\mathbb{CP}(2)$ com curva alg\'ebrica invariante e
singularidades radiais}, Instituto de Estudios con Iberoamerica y
Portugal, Seminarios Tem\'aticos, Vol {\bf 2}, Fasc. {\rm III},
Septiembre 1998.

\bibitem{LinsNeto2} {\sc Lins Neto, A.},{\it \ Some examples for the Poincar\'e and Painlev\'e
problems.}, Ann. Scient. \'Ec. Norm. Sup., $4^{e}$ s\'erie, t. {\bf 35} (2002), p.231--266.

\bibitem{lunkevichsibirskii} {\sc Lunkevich, V.A.; Sibirski\u{\i},
K.S.},{\it \ Integrals of a general quadratic differential system
in cases of a center.} Diff. Equations {\bf 18} (1982), 563--568.

\bibitem{Moulin} {\sc Moulin Ollagnier, J.},{\it \ About a
conjecture on quadratic vector fields.} J. Pure Appl. Algebra {\bf
165} (2001), 227--234.

\bibitem{Moussu} {\sc Moussu, R.}, {\it \ Une d\'emonstration d'un th\'eor\`eme de Lyapunov-Poincar\'e},
Ast\'erisque {\bf 98-99} (1982), 216--223.

\bibitem{Poincare5} {\sc Poincar\'e, H.},{\it \ Sur
l'int\'egration alg\'ebrique des \'equations diff\'erentielles du
premier ordre et du premier degr\'e,} Rend. Circ. Mat. Palermo
{\bf 5} (1891), 161--191.

\bibitem{prellesinger} {\sc Prelle, M.J.; Singer, M.F.},{\it \ Elementary
first integrals of differential equations.} Trans. Amer. Math.
Soc. {\bf 279} (1983),  215--229.

\bibitem{Schlomiuk1} {\sc Schlomiuk, D.},{\it \ Algebraic paticular integrals, integrability and the problem of center.}, Trans. Amer. Math. Soc. {\bf 338} (1993), 799-841.

\bibitem{Schlomiuk2} {\sc Schlomiuk, D.},{\it \ Algebraic and geometric
aspects of the theory of polynomial vector fields.} Bifurcations
and periodic orbits of vector fields (Montreal, PQ, 1992),
429--467, NATO Adv. Sci. Inst. Ser. C Math. Phys. Sci., {\bf 408},
Kluwer Acad. Publ., Dordrecht, 1993.

\bibitem{singer} {\sc Singer, M.F.},{\it \ Liouvillian first
integrals of differential equations.} Trans. Amer. Math. Soc. {\bf
333} (1992), 673--688.

\bibitem{wolfram} {\texttt{ http://functions.wolfram.com}}, \  2000-2003 Wolfram Research, Inc.

\bibitem{Z1} {\sc \.{Z}o\l\c{a}dek, H.},{\it \ The classification of reversible cubic systems
with center.} Topol. Methods Nonlinear Anal. {\bf 4} (1994),
79--136.

\bibitem{Z2} {\sc \.{Z}o\l\c{a}dek, H.},{\it \ Remarks on: ``The classification of reversible cubic
systems with center.''} Topol. Methods Nonlinear Anal. {\bf 8}
(1996), 335--342.

\end{thebibliography}
\end{document}